\documentclass[11pt]{article}
\pdfoutput=1
\usepackage[ascii]{inputenc}
\usepackage[T1]{fontenc}
\usepackage[english]{babel}
\usepackage{amsmath,amssymb,amsfonts,textcomp,epsfig,amsthm}
\usepackage[pdftex,colorlinks=true,citecolor=blue]{hyperref}
\usepackage{color}
\usepackage{calc}
\usepackage[all]{xy}
\usepackage{graphicx}
\usepackage{tikz}
\usepackage{appendix}

  \newcommand{\R}{\mathbb{R}}
\newcommand{\C}{\mathbb{C}}

\newcommand{\Z}{\mathbb{Z}}

\newtheorem{thm}{Theorem}

\newtheorem{coro}{Corollary}

\newtheorem{proposition}{Proposition}
\newtheorem{lem}{Lemma}

\newenvironment{pf}[1][Proof]{\begin{trivlist}
\item[\hskip \labelsep {\bfseries #1}]}{\end{trivlist}}

\title{Special cubic perturbations of the Duffing oscillator $x''=x-x^3$ near the eight-loop }
\author{
Lubomir Gavrilov\\
Institut de Math\'{e}matiques de Toulouse, UMR 5219\\
Universit\'{e}  de Toulouse,  31062 Toulouse,  France.\\
E-mail: \texttt{lubomir.gavrilov@math.univ-toulouse.fr}\\
\\
Ameni Gargouri\\
Facult\'e des Sciences de Sfax, D\'epartement de Math\'ematiques,\\
BP 1171, 3000 Sfax, Tunisie.\\
E-mail: \texttt{ameni.gargouri@gmail.com}\\
\\
Bassem Ben Hamed\\
Ecole Nationale d' Electronique et des T\'el\'ecommunications de Sfax\\
Technop\^ole El Ons, Route de Tunis km 10, BP 1163, 3021 Sfax, Tunisie.\\
E-mail: \texttt{bassem.benhamed@gmail.com}\\
}

\date{}
\begin{document}
\maketitle
\date
\begin{abstract}
We find an upper bound for the number of limit cycles, bifurcating from the 8-loop of the Duffing oscillator $x''= x-x^{3}$ under the special cubic perturbation
$$
x''= x-x^{3}+\lambda_{1}y+\lambda_{2}x^{2}+\lambda_{3}xy+\lambda_{4}x^{2}y$$.
\end{abstract}

\section{Introduction}
\label{section0}

\noindent

The perturbed Duffing oscillator 
\begin{eqnarray}\label{eq1}
X_{\lambda}:
\left\{\begin{array}{ccl} \dot{x}&=&y \\
\dot{y}&=&x-x^{3}+\lambda_{1}y+\lambda_{4}x^{2}y 
\end{array}\right.
\end{eqnarray}
where $\lambda_i$ are small parameters, appears in a
unavoidable way in the study of co-dimension two versal deformations of plane vector field with two zero eigenvalues and symmetry of order 2,
\cite[Horozov, 1979]{horo79}, \cite[Carr, 1981]{carr81}.  
For $\lambda_1=\lambda_3 = 0$ (\ref{eq1}) has a figure eight loop (union of two homoclinic saddle connections) which is the level set $\{H=0\}$ of the first integral $H$ (\ref{hamiltonian}). 
It is shown in the above mentioned papers that the Melnikov  function responsible for the bifurcations  has at most one zero near the figure eight loop, respectively  at most one limit cycle tends to the figure eight  loop of the non-perturbed system (\ref{eq1})  ( for a self-contained proof see \cite[section 4.2.]{clw94}).
The  simplicity of the  the corresponding bifurcation diagrams hides at least two more complicated phenomena, which might appear in deformations of higher co-dimension. 
\begin{itemize}
\item First, when the (first) Melnikov function vanishes, the system need not be "integrable". 
\end{itemize}
Such is the case
with the following more general deformation, 
\begin{eqnarray}\label{f4}
X_{\lambda}:
\left\{\begin{array}{ccl} \dot{x}&=&y \\
\dot{y}&=&x-x^{3}+\lambda_{1}y+\lambda_{2}x^{2}+\lambda_{3}xy+\lambda_{4}x^{2}y 
\end{array}\right.
\end{eqnarray}
where $\lambda_i$ are small parameters. As we shall see, the Bautin ideal of the first return map related to the exterior period annulus is generated by
 $\lambda_1, \lambda_4$ and $\lambda_2 \lambda_3$ (Theorem \ref{bautin}). This implies that if $\lambda_1=\lambda_4=0$ but $\lambda_2\lambda_3\neq 0$, the first Melnikov function is identically zero, but the exterior period annulus is destroyed. The higher Melnikov functions were  computed by Iliev (1998) \cite{Iliev98}, and it follows from their analysis that at most two zeros   can bifurcate near the 8-loop. 
 \begin{itemize}
\item Second, not all limit cycles near the 8 loop need to be "shadowed" by zeros of Melnikov function.
\end{itemize}
The existence of such  "alien" limit cycles was  discovered by Dumortier, Roussarie and Caubergh (2005)  \cite{F. Dum, M. Cau}.
Up to these papers, the one-to-one correspondance between zeros of the Melnikov function and limit cycles near polycycle more complicated than a homoclinic loop were sometime considered as granted, leaving some place for speculations.

In the present paper we bound the cyclicity of the 8-loop with respect to the deformation (\ref{f4}). In contrast to (\ref{eq1})we study here a co-dimension three deformation. The "expected" number of limit cycles bifurcating from the 8-loop is therefore two. Indeed, when the first Melnikov function $M_1\neq 0$ we show that it is indeed the case. But when $M_1=0$ we show that this cyclicity is at most five. This bound is presumably bigger than the cyclicity of the 8-loop, but yet it is the only available bound.

\begin{figure}
\begin{center}
 \def\svgwidth{12cm}
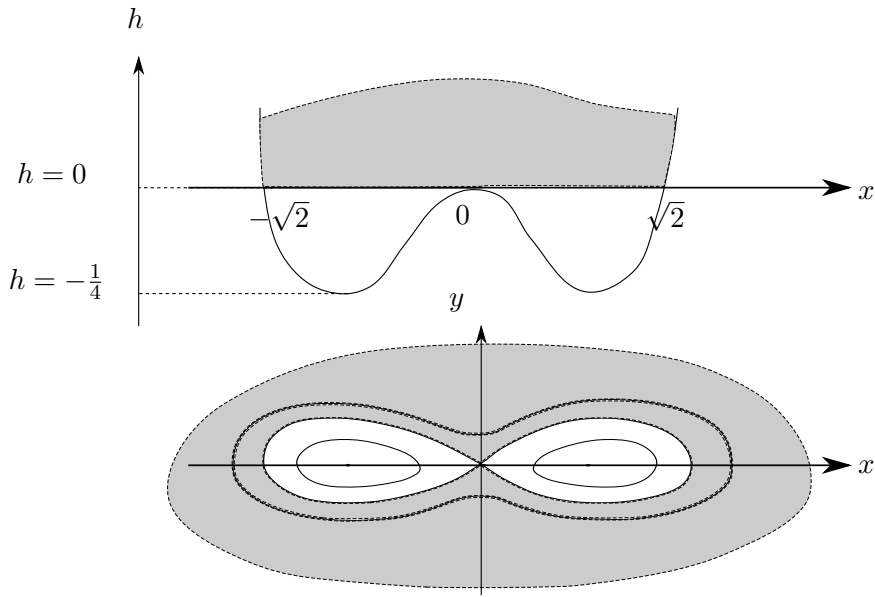
\end{center}
\caption{Phase portrait of $X_0$ on the $(x,y)$-plane and the graph of  $H(x,0)=- \frac{x^2}{2} + \frac{x^4}{4}$}
\label{fig1}
\end{figure}

The system (\ref{f4})  has been studied by many authors, see \cite{Iliev98,Iliev99,lmr06,clw94}
and the references given there.
The phase portrait of $X_0$, which has a first integral 
\begin{equation}
\label{hamiltonian}
H(x,y) = \frac{y^2}{2} - \frac{x^2}{2} + \frac{x^4}{4}
\end{equation}
is shown on Fig.\ref{fig1}. As a first approach one may consider one-parameter deformations of $X_0$.
More precisely,  consider  analytic arcs
$$
\varepsilon \mapsto \lambda(\varepsilon), \lambda(0)=0
$$
in the parameter space $\{\lambda_i \}$, and the corresponding one-parameter deformation $X_{\lambda(\varepsilon)}$.

To each annulus of $X_0$ and   deformation $X_{\lambda(\varepsilon)}$    one may associate the
Poincar\'e return map
\begin{equation}
\label{returnmap}
P_{\varepsilon}(h)=h+\varepsilon^{k}M_{k}(h)+O(\varepsilon^{k+1}), M_{k}\neq0
\end{equation}
where the zeros of the k-th order Melnikov function $M_k$ control the limit cycles of $X_\varepsilon$. 

Zoladek and Jebrane and 
Iliev and Perko \cite{Iliev99} computed the Melnikov functions $M_k$ and studied their zeros. Based on this they found 
the cyclicity of the three open period annuli (two interior and one exterior annulus) of $X_0$ with respect to the one-parameter deformation
$X_{\lambda(\varepsilon)}$. Later Li, Mardesic and Roussarie \cite{lmr06} extended  these results to multi-parameter deformations.
Indeed, according to \cite[Theorem 1]{gavr08}, the study of one-parameter deformations $X_{\lambda(\varepsilon)}$ is enough to estimate the cyclicity of the open annuli of the multi-parameter deformation $X_\lambda$. This is a general result,  close to the way in which  the Nash space of arcs of a singular variety is used to desingularize it, see \cite{gafr} for details.

The results in the above mentioned papers provide estimates to the cyclicity not only of the period annuli, but also 
 of the two homoclinic loops of $X_0$, as it  follows from a classical result of Roussarie \cite{rous98}. 
 
 The present paper is   devoted to the study of the missing cyclicity of the 8-loop, see
Fig. \ref{ii}, e.g. \cite[section 5]{lmr06}. 

We use complex methods, in the spirit of \cite{e14,e15,Gav13}. Our main result is that at most five limit cycles can bifurcate from eight-loop  (Theorem \ref{main}), although we did not succeed to prove that this bound is exact. It is interesting to note, that even for a generic perturbation \eqref{f4}, two limit cycle can appear near a eight-loop, while at the same time the first Melikov function exhibits only one zero.
Hence there is a limit cycle that is not covered by a zero of the related Abelian integral. Such a limit cycle were called  "alien" in  \cite{F. Dum, M. Cau}. Some partial results in this sense can be found in \cite{yyh14}. 

Instead of polynomial, one may consider general analytic families of analytic vector fields $X_\lambda$, such that $X_0$ has an eight loop. The finite cyclicity of the eight loop in this context follows from \cite{e15}. The main difficulty to prove this result is the case, when the return map $P_0$ associated to the eight-loop of $X_0$ is the identity map (eight-loop of infinite co-dimension). If the return map is not the identity map (the case of finite co-dimension) the finite cyclicity  together with an explicit bound for the number of limit cycles has been found by Jebrane and Mourtada \cite{mourtada}. Their result remains true in the $C^\infty$ category.


The paper is organized as follows. In the next section \ref{rappelarticle1} we describe the Bautin ideal associated to (\ref{f4}) and its \emph{exterior period annulus}, as well the corresponding 
Melnikov functions. These results are classical. In section \ref{sectionRES} we formulate our main result - Theorem \ref{main}, which says that the cyclicity of the eight-loop is at most five.
Its proof is based on Theorem \ref{th81} in which we study one-parameter deformations of $X_0$. The proof of Theorem \ref{th81} is carried out in section \ref{sectionth4}.

\section{The Bautin ideal and the Melnikov functions}
\label{rappelarticle1}

Define  the complete elliptic integrals
$$
I_{i}(h) = \oint_{\gamma(h)}x^{i}ydx, I'_{i}(h) = \oint_{\gamma(h)}\frac{x^{i}}{y}dx.
$$
where
\begin{equation}
\label{gamma}
\gamma(h) = \{ (x,y)\in \R : H(x,y) = h > 0 \}.
\end{equation}

Then the first return map $P_{\varepsilon}$, see (\ref{returnmap}), near an oval $\gamma(h)$ is well defined and for the 
 first non-vanishing Poincar\'e-Pontryagin-Melnikov function $M_{k}$  we have
\begin{thm}[Iliev, \cite{Iliev98}]  
\label{melnikov}
 If  $M_1(h)$ is not identically zero, then
\begin{equation}
\label{first}
M_1(h) = \lambda_1 I_0(h) + \lambda_4 I_2(h), \lambda_i \in \R
\end{equation}
otherwise  
\begin{equation}
\label{principal}
M_k(h)= \lambda_{1k}I_0(h)+ \lambda_{4k}I_2(h) + \frac13 \sum_{i+j=k} \lambda_{2i}\lambda_{3j} I_4'(h),\;\; k\geq 2 .
\end{equation}
\end{thm}
\label{center}

Let $\delta(h) \subset \{ (x,y)\in \C^2 : H(x,y)=h\}$ be a continuous family of closed loops, representing a cycle in $H_1(\{H(x,y)=h\}, \Z)$ which vanishe at the saddle point when $h$ tends to $0$. Note that $\delta(h)$ can not be represented in a real domain. It is the variation of the real oval
 $\gamma(h) \in H_1(\{H(x,y)=h\}, \Z)$ when $h$ makes one turn around $0\in \C$ in a complex domain.  Define 
$$
\tilde I_{i}(h) = \oint_{\delta(h)}x^{i}ydx, \tilde I'_{i}(h) = \oint_{\delta(h)}\frac{x^{i}}{y}dx.
$$

 Similarly to $P_\varepsilon$, define $\tilde{P}_\varepsilon$ to be the holonomy map of one of the separatrices of the saddle point and write
 \begin{equation}
\label{holonomy}
\tilde{P}_{\varepsilon}(h)=h+\varepsilon^{k}\tilde{M}_{k}(h)+O(\varepsilon^{k+1}), \tilde{M}_{k}\neq 0 .
\end{equation}

Repeating the proof of Theorem \ref{melnikov} we obtain (see also \cite{Iliev99})
\begin{thm}
\label{melnikov1}
 If   $\tilde M_1(h)$ is not identically zero, then
\begin{equation}
\label{firstt}
\tilde{M}_1(h) = \lambda_1 \tilde{I}_0(h) + \lambda_4 \tilde{I}_2(h), \lambda_i \in \R
\end{equation}
otherwise  
\begin{equation}
\label{principalt}
\tilde M_k(h)= \lambda_{1k}\tilde I_0(h)+ \lambda_{4k}\tilde I_2(h) + \frac13 \sum_{i+j=k} \lambda_{2i}\lambda_{3j} \tilde I_4'(h),\;\; k\geq 2 .
\end{equation}
\end{thm}

The next two Lemmas are proved in a standard way and can be found in \cite{Iliev99,p}.
\begin{lem}
\label{lemmarea}
\noindent\\The complete elliptic integrals integrals $I_0, I_1, I_2$, satisfy the following Picard-Fuchs system :
\begin{eqnarray}\label{area}
I_0(h)&=&\frac{4}{3}hI'_0(h)+\frac{1}{3}I'_2(h) \\
I_2(h)&=&\frac{4}{15}hI'_0(h)+\left(\frac{4}{5}h+\frac{4}{15}\right)I'_2(h)\\
(4h+1)I'_4(h)&=&4hI_0(h)+5I_2(h)\\
4h(4h+1)I''_0(h)&=&-3I_0(h).
\end{eqnarray}
\end{lem}
Based on Lemma \ref{lemmarea} it can be deduced
\begin{lem}\label{14}
\noindent  The complete elliptic integrals integrals $I_0, I_1, I_2$ allow   convergent   expansions near $h=0$ of the form
\begin{eqnarray*}
I_{0}(h)&=&(-h+\frac{3}{8}h^{2}-\frac{35}{64}h^{3}+...)\ln h+\frac{4}{3}+a_{1}h+a_{2}h^{2}+...\\
I_{2}(h)&=&(\frac{1}{2}h^{2}-\frac{5}{8}h^{3}-\frac{315}{256}h^{4}...)\ln h+\frac{16}{15}+4h+b_{2}h^{2}+...\\
I'_{4}(h)&=&(-\frac{3}{2}h^{2}+\frac{35}{8}h^{3}-\frac{471}{256}h^{4}+...)\ln h+\frac{16}{3}+4h+(4a_{1}+5b_{2}-\frac{304}{3})h^{2}+...\\
\end{eqnarray*}
where $a_1, a_2, b_2 $ are constants.
\end{lem}
\begin{coro}
\label{cor1}
The first non-vanishing Melnikov function $M_k$, $k\geq 1$, allows a convergent  expansion of the form
\begin{equation}\label{MK1}
M_{k}(h)=c_{0}+c_{1}h\ln h+c_{2}h+c_{3}h^{2}\ln h +...
\end{equation}
Moreover, if $k=1$ and $c_0=0$ then $c_1\neq 0$.
If $k\geq 2$ and $c_0=c_1=0$ then $c_2 \neq 0$. 
\end{coro}
Lemma \ref{14}  implies  
\begin{lem}\label{141}
\noindent  The complete elliptic integrals integrals $\tilde I_0, \tilde I_1, \tilde I_2$ are analytic near $h=0$ and
allow       expansions   of the form
\begin{eqnarray*}
 \tilde I_{0}(h)&=&2\pi \sqrt{-1}  (-h+\frac{3}{8}h^{2}-\frac{35}{64}h^{3}+...) \\
\tilde I_{2}(h)&=&2\pi \sqrt{-1}  (\frac{1}{2}h^{2}-\frac{5}{8}h^{3}-\frac{315}{256}h^{4}...) \\
\tilde I'_{4}(h)&=&2\pi \sqrt{-1}  (-\frac{3}{2}h^{2}+\frac{35}{8}h^{3}-\frac{471}{256}h^{4}+...)  .
\end{eqnarray*}
\end{lem}

To define the Bautin ideal, we note that the first return map associated to the exterior annulus of $X_0$ is also defined for all $\lambda$ close to $0$. Take any $h_0>0$ and expand the Poincar\'e return map
$$
P_\lambda(h) - h = \sum_{k\geq 0} \beta_k (\lambda) (h-h_0)^k
$$
where $\beta_k= \beta_k(\lambda)$ are suitable analytic functions which generate an ideal $B(h_0)$ 
 in the Noetherian ring of germes of analytic functions $\mathcal O(\R^4,0)$. It is known that the ideal $B(h_0)$ does not depend on $h_0>0$, and it is called the Bautin ideal associated to the deformed vector field $X_\lambda$ and to the exterior period annulus.  
 \begin{thm}
 \label{bautin}
 The Bautin ideal $\mathbb B$ is  generated by the polynomials $\lambda_1, \lambda_4$ and $\lambda_2 \lambda_3$.
 \end{thm}
 The proof of this remarkable fact follows from Theorem \ref{melnikov}, and can be also found in \cite{lmr06}. Indeed, the Poincar\'e return map $P_\varepsilon$ 
 is the identity map, if and only if $M_k=0, \forall k \geq 1$. This implies that
if the deformed vector field $X_{\lambda(\varepsilon)} $ 
allows for all small $\varepsilon$ a continuous band of periodic orbits on the exterior of a 8-loop if and only if
\begin{equation}
\label{lambdaeps}
\lambda_1(\varepsilon) = \lambda_4(\varepsilon) = \lambda_2(\varepsilon) \lambda_3(\varepsilon) = 0 .
\end{equation}
The center variety defined by the Bautin ideal $\mathbb B$ is a germ of analytic set centred at the origin
 $0\in \C^4$, which combined with (\ref{lambdaeps}) implies
$$
C = \{ \lambda : \lambda_1=\lambda_4= \lambda_2 = 0 \}\cup \{ \lambda : \lambda_1=\lambda_4= \lambda_3 = 0 \}
$$
 This already shows that $\mathbb B$ is polynomially generated, although it  does not need to be radical. As the ideal of the center variety $C$ is
generated by $\lambda_1, \lambda_4, \lambda_2 \lambda_3$ then we can divide the displacement map $P_\lambda(h) - h$ and write
$$
P_\lambda(h) - h = \lambda_1 f_1 + \lambda_4 f_2 +  \lambda_2 \lambda_3 f_3
$$
where $f_i = f_i(\lambda, h-h_0)$ are germs of analytic functions. Substituting $\lambda=\lambda(\varepsilon)$ and taking into consideration
Theorem \ref{melnikov} we conclude that 
\begin{equation}
\label{divide}
P_\lambda(h) - h = \lambda_1 (I_0(h) + O(\lambda)) +
 \lambda_4( I_2(h)+O(\lambda)) + \frac13  \lambda_2 \lambda_3 (I_4'(h) + O(\lambda) )
\end{equation}
where, by abuse of notations, $O(\lambda)$ is a germ of analytic function in $\lambda$ and $h-h_0$ which vanishes at $\lambda=0$.

\section{Cyclicity of the 8-loop}
\label{sectionRES}
The so called 8-loop is the union of the two homoclinic orbits of the Hamiltonian system $X_0$ having a first integral $H(x,y)$ as on Fig.\ref{fig1}.
The cyclicity of the 8-loop is the maximal number of limit cycles of $X_\lambda$, which tend to to the 8-loop, when $\lambda= (\lambda_1, \lambda_2,\lambda_3,\lambda_4)$ tends to $0$, for rigorous definition see Roussarie \cite{rous98}.
The main result of the paper is 
\begin{thm}
\label{main}
The cyclicity of the 8-loop of $X_0$ with respect to the four-parameter family of cubic deformations $X_\lambda$, defined by (\ref{f4}), is at most equal to five.
\end{thm}

Following \cite{Iliev99}, instead of $\lambda \in \R^4$ we consider first analytic arcs (one-parameter deformations)
$$
\varepsilon \mapsto \lambda(\varepsilon) = (\sum_{k\geq 1} \lambda_{1k}\varepsilon ^k,  \sum_{k\geq 1} \lambda_{2k}\varepsilon ^k, \sum_{k\geq 1} \lambda_{3k}\varepsilon ^k, \sum_{k\geq 1} \lambda_{4k}\varepsilon ^k) .
$$
together with the corresponding one-parameter family of vector fields 
$$\varepsilon \mapsto X_{\lambda(\varepsilon)} .$$
To the exterior period annulus of $X_0$ we associate a Poincar\'e first return map (the identity map) which is well defined also for $\lambda$ close to $0$. As usual we parameterize this map by the restriction $h$ of $H(x,y)$ on a cross section to the exterior period annulus of $X_0$ and write
$$
P_{\lambda(\varepsilon)} (h) = h+\varepsilon^{k}M_{k}(h)+O(\varepsilon^{k+1}), M_{k}\neq 0 .
$$

For such one-parameter deformation we can give a more precise result
\begin{thm}
\label{th81}
 The cyclicity of the 8-loop of $X_0$ with respect to the one-parameter deformation $X_{\lambda(\varepsilon)}$ is at most equal to 
  \begin{itemize}
\item 
two, if $M_{1}\neq 0$  
\item five, otherwise .
\end{itemize}
\end{thm}

\begin{figure}
\begin{center}
 \def\svgwidth{6cm}
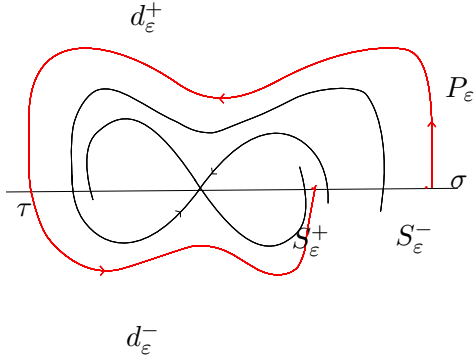
\end{center}
\caption{Monodromic eight-loop and the Dulac map $d^{\pm}_{\varepsilon}$}
\label{ii}
\end{figure}

\subsection{Proof of Theorem \ref{th81}}
\label{sectionth4}

We adapt the proof given for the so called two-saddle loop in \cite{e14}, where the reader will find more detailed theoretical justification of the method which we use. We outline first the plan of our proof.

\begin{figure}
\begin{center}
\includegraphics{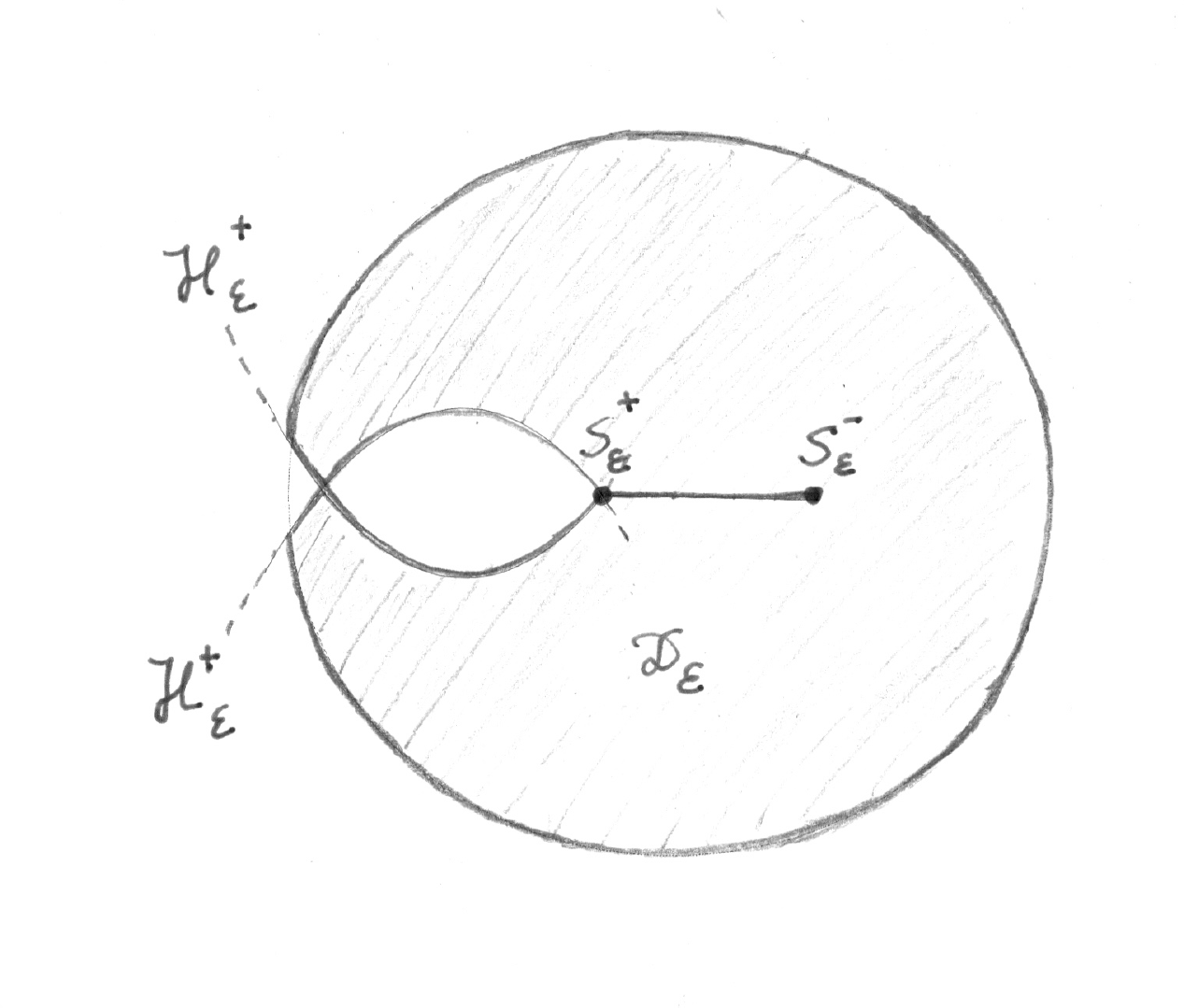}
\end{center} 
\caption{The domain $\mathcal{D}_\varepsilon$}
\label{cde}
\end{figure}

\noindent Consider the Dulac maps $d^{+}_{\varepsilon}$, $d^{-}_{\varepsilon}$ associated to the perturbed foliation,
and to the cross sections $\sigma$ and $\tau$, see Fig. \ref{ii}. We parameterize each cross-section by the
restriction of the first integral \emph{f} on it, and denote $h= \emph{f}\mid \sigma$. Each function $d^{\pm}_{\varepsilon}$
is multivalued and has a critical point at $S^{\pm}_{\varepsilon}\in \R$,  $S^{\pm}(0)= 0$. The  points $S^{+}$, $S^{-}$
depend analytically on $\varepsilon$. Without loss of generality we shall suppose that
$\varepsilon> 0 $ and $S^{-}_{\varepsilon} > S^{+}_{\varepsilon}$, see Fig. \ref{cde}. A limit cycle intersects the cross-section  $\sigma$ at $h$
if and only if  $d^{+}_{\varepsilon}(h)= d^{-}_{\varepsilon}(h)$.  The Poincar\'e  return map $P_{\varepsilon} = P_{\lambda(\varepsilon})$ is defined as
\begin{eqnarray}\label{pf}
 P_{\varepsilon}= d^{+}_{\varepsilon}\circ (d^{-}_{\varepsilon})^{-1} .
\end{eqnarray}
The limit cycles of $X_{\varepsilon}= X_{\lambda(\varepsilon)}$ correspond also to the fixed points of $P_{\varepsilon}$, the zeros of 
$$P_{\varepsilon} - id =  (d^{+}_{\varepsilon}-d^{-}_{\varepsilon})\circ (d^{-}_{\varepsilon})^{-1}$$
as well the  zeros of 
the dispalacment map $d^{+}_{\varepsilon}-d^{-}_{\varepsilon}$. Each of the Dulac maps $d^{+}_{\varepsilon}, d^{-}_{\varepsilon}$ has a single singular point $S^\pm_\varepsilon$ corresponding to the saddle point of $X_\varepsilon$
and otherwise  allows an analytic continuation in a complex domain to a multivalued function. Assume that $S^+_\varepsilon < S^-_\varepsilon$ as on Fig. \ref{cde}. To count the zeros of the displacement map on the interval $(S^-_\varepsilon, R)$ we shall bound them by the number of the zeros of the displacement map in the larger complex domain $\mathcal D_\varepsilon$ which is shown on Fig. \ref{cde}.  We recall from \cite{e14} that the domain $\mathcal D_\varepsilon$ is bounded by a small circle $S_R$ of radius $R$, by the segment 
 $(S^{+}_{\varepsilon}, S^{-}_{\varepsilon})$, and by the zero locus $\mathcal H^+_\varepsilon$ of the
imaginary part of $d^{+}_{\varepsilon}$
$$
\mathcal H^+_\varepsilon = \{z\in \C: \operatorname{Im} d^{+}_{\varepsilon} (z) = 0 \} .
$$
Define similarly
$$
\mathcal H^-_\varepsilon = \{z\in \C: \operatorname{Im} d^{-}_{\varepsilon} (z) = 0 \} 
$$
and recall from \cite{e14}  that $\mathcal H^+_\varepsilon$, $\mathcal H^-_\varepsilon$ as subsets of $\C=\R^2$ are (germs of) real analytic curves.

The number of the zeros of $d^{+}_{\varepsilon}-d^{-}_{\varepsilon}$ in ${\cal D_{\varepsilon}}$ is computed according
to the argument principle: it equals the increase of the argument along the boundary of ${\cal D_{\varepsilon}}$.

\noindent $\bullet$ Along the circle and far from the critical points, the displacement function is
"well" approximated by $\varepsilon^{k}M_{k}(h)$ which allows one to estimate the increase of the
argument.\\
\noindent $\bullet$ Along the segment $(S^{+}_{\varepsilon}, S^{-}_{\varepsilon})$ the zeros of the imaginary part of the displacement function coincide with the fixed points of the holomorphic holonomy map along
the separatrix through $S^{-}_{\varepsilon}$. The zeros are therefore well approximated,  by the  Abelian integral $\tilde M_k$, along the vanishing cycle $\delta(h)$.

\noindent $\bullet$ 
Along the zero locus of the imaginary part of $d^{+}_{\varepsilon}$, the zeros of the imaginary part of the displacement map $d^{+}_{\varepsilon} - d^{-}_{\varepsilon}$ coincide with $
\mathcal H^+_\varepsilon  \cap \mathcal H^-_\varepsilon$, which are in fact the fixed points of suitable holonomy map, which
can be desribed by analogy to  \cite{e15}.

\begin{figure}[htbp]
\begin{center}
\includegraphics[width=12cm]{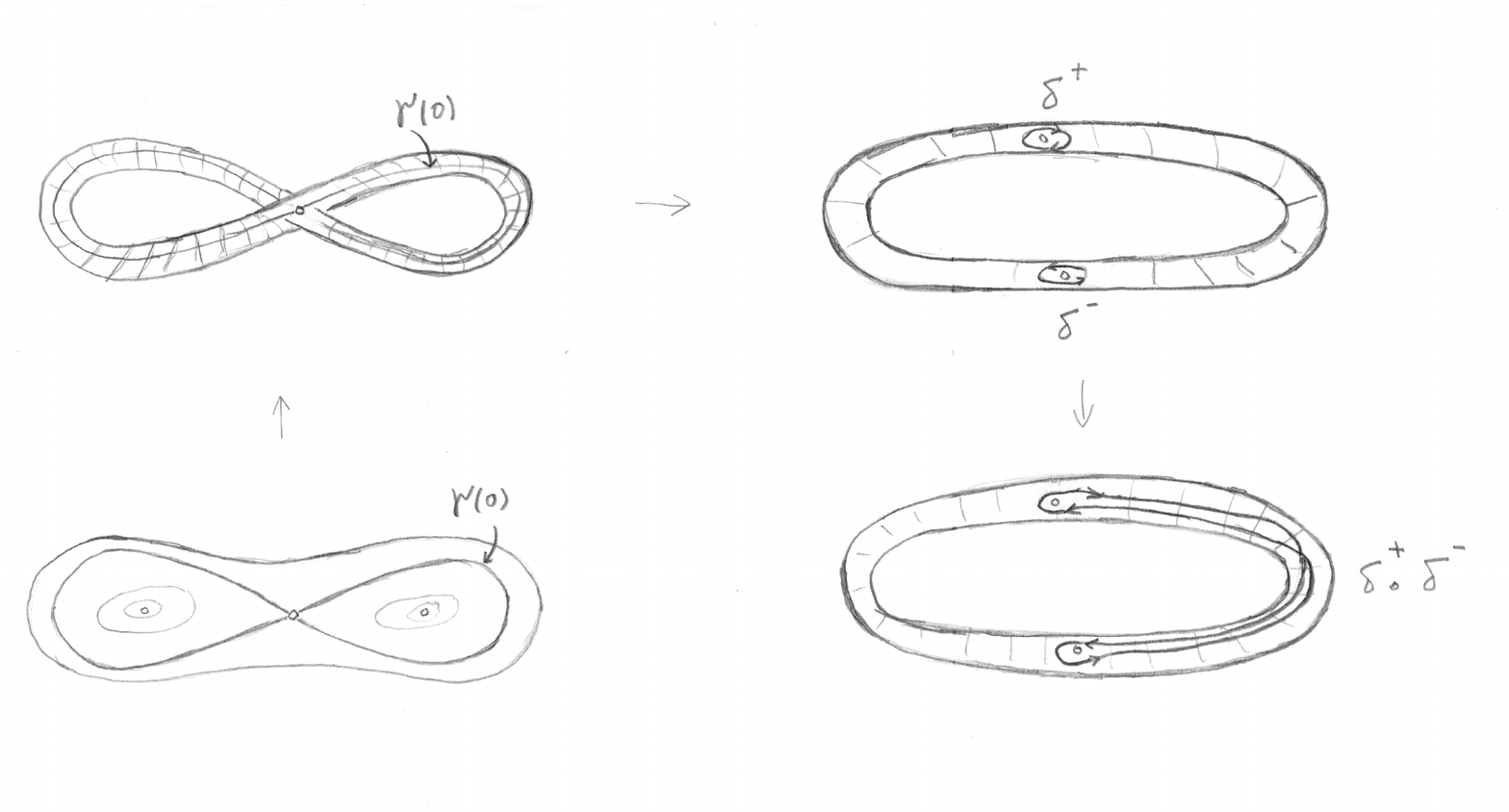}
\caption{We show successively the eight-loop $\gamma(0)$ of the vector field $X_0$, the complex neighboourhood  $U \subset \{ H(x,y)= 0 \}$  of  $\gamma(0)$, the two closed loops $\delta^\pm$ in the neighbourhood $U$, the composed loop $\delta^+\circ \delta^-$ in $U$.}
\label{a8}
\end{center}
\end{figure}

The last point needs some explication. Namely, to a closed loop $l$ contained in a leaf of $X_0$ we associate a holonomy map $\textbf{h}^{\varepsilon}_{l}$ of the perturbed system $X_\varepsilon$. 
Let $\delta^+, \delta^-$ be two closed loops in the two local separatrices of $X_0$ through the saddle point $(0,0)$, that is to say
$\delta^\pm \subset \{ (x,y) \in \C^2 : H(x,y) = 0 \}$ and let $U \subset \{ H(x,y)= 0 \}$ be a complex neighboorhood of the real eight-loop, which topologically is an
annulus with two identified marked points $S^\pm$ identified to a single point, which is the saddle $(0,0)$, shown on Fig. \ref{fig1}. The loops $\delta^\pm$ are considered up to a free homotopy and denote $\delta^+ \circ \delta^-$ the composed loop as it is
 shown on Fig. \ref{a8}. Note that the loops $\delta^\pm$ are contractible in $U$, and that they allow a continuation to a family of loops $\delta(h)^+,\delta^-(h)$ which are homologous on $H_1(\{H=h\},\Z)$. Therefore they define to a continuous family of vanishing cycles $\delta(h)$.

We have defined in such a way three holonomy maps 
$$
\textbf{h}^{\varepsilon}_{\delta ^{-}}, \textbf{h}^{\varepsilon}_{\delta ^{+}}, \textbf{h}^{\varepsilon}_{\delta ^{-}\circ \delta ^{+}} =  \textbf{h}^{\varepsilon}_{\delta ^{+}}\circ \textbf{h}^{\varepsilon}_{\delta ^{-}} .
$$
The fixed points of $\textbf{h}^{\varepsilon}_{\delta ^{+}}\circ \textbf{h}^{\varepsilon}_{\delta ^{-}}$ coincide with $\mathcal H^+_\varepsilon  \cap \mathcal H^-_\varepsilon$, that is to say with the zeros of the imaginary part of $d^{+}_{\varepsilon} - d^{-}_{\varepsilon}$
along the the imaginary part of $d^{+}_{\varepsilon}$.

To the end of the section we follow the steps outlined above, by completing the missing estimates.

\subsubsection{The case $M_{1}\neq0$}
\label{sectionecase2}
\noindent In this section we consider the perturbed   system (\ref{f4}) under the generic assumption that
$$M_{1}(h)=\int_{\gamma(h)}\omega|_{\varepsilon=0}=\lambda_1 I_0(h) + \lambda_4 I_2(h), \lambda_i \in \R$$
is not identically zero.\\
\noindent By Lemma \ref{14}  the Abelian integrals $I_1, I_4$ are linearly independent and hence $M_1\neq 0$ if and only if $(\lambda_1,\lambda_4)\neq (0,0)$.
 The Poincar\'e-Pontryagin function $M_{1}$ has a continuous limit  at $h=0$ to $M_1(0)$ whichis the classical Melnikov integral along the eight-loop $\gamma(0)$. 
 It is  known that  the vanishing of the Melnikov integral   $M_{1}(0)$ is a necessary condition for a bifurcation of a limit cycle :

\begin{proposition}\label{proppp1}
If $M_{1}(0)\neq0$, then no limit cycles bifurcate from the eight-loop $\gamma$.
\end{proposition}

\begin{pf}
Suppose that there is a sequence of limit cycles  $ \{\delta_{\varepsilon_{i}}\}_i$  of $X_{\lambda(\varepsilon_{i})} $ which tend to the eigth-loop  $\gamma(0)$ and  $\varepsilon_{i}\to 0$ when $i\to \infty$. Then
$$0=-\int_{\delta_{\varepsilon_{i}}}dH=\varepsilon_{i}\int_{\delta_{\varepsilon_{i}}}\omega$$
which implies
$$0=lim_{\varepsilon_{i}\rightarrow0}\int_{\delta_{\varepsilon_{i}}}\omega=\int_{\gamma(0)}\omega|_{\varepsilon=0}=M_{1}(0)$$
$\Box$. 
\end{pf}

By Corollary \ref{cor1} at most one zero of $M_1(h)$ bifurcates from $h=0$. It can be proved, however, that two limit cycles can bifurcate from the eight loop, when $\lambda_1, \lambda_4$ tend to zero. Thus, an alien limit cycle is present near the eight-loop, see 
 \cite{F. Dum, M. Cau, Gav13}.
We shall prove here the following weaker 
\begin{proposition}\label{prope1}
If the first Melnikov function is not identically zero, then at most two limit cycles bifurcate from the eight-loop $\gamma(0)$.
\end{proposition}
\begin{itemize}
\item
By Proposition \ref{proppp1}, if limit cycles bifurcate from the eight-loop, then
$c_{0}= \lambda_{1}I_{0}(0)+\lambda_{4}I_{2}(0)=0$
and hence 
$ \lambda_{1} +4 \lambda_{4} = 0$. As
$$
(d^{+}_{\varepsilon}-d^{-}_{\varepsilon})(h)=(P_{\varepsilon}-id) \circ d^{-}_{\varepsilon}(h) = \varepsilon M_{1}(h) + \dots
$$
then the displacement map along the circle $S_{R}$ is approximated by $\varepsilon M_{1}$ which has as a leading term $h \ln h$
(because if $c_{0}=0$ then $c_{1} \neq 0$). The increase of the argument of $h\ln h$, and hence
of the displacement map, along the circle $S_{R}$ is close to $2\pi$ but strictly less than
$2\pi$.
\vskip0.2cm

\item The imaginary part of the displacement map, along the interval $[S^{+}({\varepsilon}),  S^{-}({\varepsilon})]$
equals the imaginary part of $d^{-}_{\varepsilon}(h)$. Its zeros equal the number of intersection
points of $\emph{H}^{+}_{\varepsilon}$ with the real axes, that is to say the fixed points of the holonomy map $\textbf{h}^{\varepsilon}_{\delta ^{-}}$, where
$$
\textbf{h}^{\varepsilon}_{\delta ^{-}} (h) = h + \varepsilon \tilde M_1(h) + \dots
$$
and $\tilde M_1(h)$ is an Abelian integral along the vanishing cycle $\delta(h)$, see Theorem \ref{melnikov1}.
The Abelian integral $M_1(h)$ has a simple zero at $h=0$ which is also a fixed point of the holonomy map.
Therefore the imaginary part of the displacement map does not vanish along the
open interval $(S^{+}({\varepsilon}),  S^{-}({\varepsilon}))$.
\vskip0.2cm
\item  The number of the zeros of the imaginary part of the displacement map, along
the real analytic curve $\mathcal H^{+}_{\varepsilon}$ equals the number of the zeros of the imaginary part
of $d^{-}_{\varepsilon}(h)$ along this curve, that is to say the number of intersection points of this curve with
$\emph{H}^{-}_{\varepsilon}$, which are the fixed points of the holonomy map $\textbf{h}^{\varepsilon}_{\delta ^{-}}\circ \textbf{h}^{\varepsilon}_{\delta ^{+}}$.

 We have 
\vskip0.2cm
\begin{equation}
\label{eqq1}
\textbf{h}^{\varepsilon}_{\delta ^{+}}(h)= h+\varepsilon \tilde M_1 (h)+O(\varepsilon^{2})
\end{equation}
\\
\begin{equation}
\label{eqq2}
\textbf{h}^{\varepsilon}_{\delta ^{-}}(h)=  h+\varepsilon \tilde M_1 (h)+O(\varepsilon^{2})
\end{equation}
\\
\begin{equation}
\label{eqq3}
\textbf{h}^{\varepsilon}_{\delta ^{+}}\circ \textbf{h}^{\varepsilon}_{\delta ^{-}}(h)= h+\varepsilon ( \tilde M_1 (h) +  \tilde M_1 (h) )+O(\varepsilon^{k+1})
\end{equation}
\\
\begin{equation}
\label{eqq4}
\textbf{h}^{\varepsilon}_{\delta ^{-}}\circ \textbf{h}^{\varepsilon}_{\delta ^{+}}(h)= h+\varepsilon ( \tilde M_1 (h) +  \tilde M_1 (h) )+O(\varepsilon^{k+1})
\end{equation}
As $ \lambda_{1} +4 \lambda_{4} = 0$, then
$\tilde M_1$ has a simple zero at $h=0$ and we conclude that the imaginary part of the displacement map vanishes at most once.
\end{itemize}
We conclude that the displacement map can have at most two zeros in the domain ${\cal D_{\varepsilon}}$ which completes the proof of Proposition \ref{prope1}.
 $\Box$

\subsubsection{The case $M_{1}=0$}
\label{sousectione3}
In this section we suppose that the Melnikov function $M_{1}(h)$ vanishes identically. The first return map has the form (\ref{returnmap}), where the Melnikov function $M_k$ is 
computed in Theorem \ref{melnikov}.
\begin{proposition}\label{prope2}
If $\lambda_{1k}\neq0$,  then at most two limit cycles bifurcate from $\gamma$.
\end{proposition}
\begin{pf}
\noindent Following the method of the preceding subsection, we evaluate the number of the
zeros of the displacement map  in the domain ${\cal D_{\varepsilon}}$
\begin{eqnarray}
d^{+}_{\varepsilon}-d^{-}_{\varepsilon}=(P_{\varepsilon}-id) \circ d^{-}_{\varepsilon} = \varepsilon^{k}M_{k}(h)+\varepsilon^{k+1}M_{k+1}(h)+...
\end{eqnarray}

\begin{itemize}
\item  Along the circle $S_{R}$
The displacement map is approximated by $\varepsilon^{k} M_{k}$ which has as a leading term: $h\ln h$ as $\lambda_{0k}\neq0$, see Lemma \ref{14}. The increase of the argument
of $h\ln h$, and hence of the displacement map, along the circle $S_{R}$ is close to $2\pi$ but strictly less than
$2\pi$.\\

\item The imaginary part of the displacement map, along the interval $[S^{+}_{\varepsilon},  S^{-}_{\varepsilon}]$
equals the imaginary part of $d^{-}_{\varepsilon}(h)$.  The number of its zeros is bounded 
 by the multiplicity of the  
Abelian integral $\tilde M_{k}(h)=\int_{\delta(h)} \omega_{k}$ having at most a simple zero at the origin, see Theorem \ref{melnikov1}.
Note, however, that the holonomy map $\textbf{h}^{\varepsilon}_{\delta ^{-}}$ has $S^{-}({\varepsilon})$ as a fixed point, and hence the cyclicity
of the saddle point is zero.
We conclude that the imaginary part of the displacement map does not vanish along the
 interval $[S^{+}_{\varepsilon},  S^{-}_{\varepsilon}]$.

\item The number of the zeros of the imaginary part of the displacement map, along
the real analytic curve $\mathcal H^{+}_{\varepsilon}$ is bounded by the the cyclicity of the zero $h=0$ of the Abelian integral
$\tilde M_k$ (see the proof of Proposition \ref{prope1}) and
$
\tilde M_k (h) = c_{1} h +  \dots , c_1 \neq 0.
$

\end{itemize}
\vskip0.2cm
Summing up the above information we conclude that at most two limit cycles bifurcate from the eight-loop.
 $\Box$
\end{pf}
\begin{proposition}\label{prope3}
If $\lambda_{1k}=0$
then at most five limit cycles bifurcate from $\gamma$.
\end{proposition}
Assuming that $M_k(0)=0$ (otherwise no limit cycles bifurcate from the eight-loop), implies that up to multiplication by a non zero constant we have
\begin{align*}
M_k(h)&= 5 I_2(h) - I_4'(h) = 16h + 4 h^2 \ln h + \dots \\
\tilde M_k(h)& = 5 \tilde I_2(h) - \tilde I_4'(h) = 8 \pi \sqrt{-1}h^2 + \dots
\end{align*}
We repeat the three steps above.
\begin{itemize}
\item  
Along the circle $S_{R}$, the displacement map is approximated by $\varepsilon^{k} M_{k}$ and hence
the increase of the argument
 of the displacement map, along the circle $S_{R}$ is close to $2\pi$.

\item 

The number of the zeros of the imaginary part of the displacement map, along the interval $[S^{+}_{\varepsilon},  S^{-}_{\varepsilon}]$
  is bounded 
 by the multiplicity of the  
Abelian integral $\tilde M_{k}(h)$ having a double zero at the origin.
Note, however, that the holonomy map $\textbf{h}^{\varepsilon}_{\delta ^{-}}$ has $S^{-}({\varepsilon})$ as a fixed point.
We conclude that the imaginary part of the displacement map   vanishes at most once along the
 interval $[S^{+}_{\varepsilon},  S^{-}_{\varepsilon}]$.

\item 
The number of the zeros of the imaginary part of the displacement map, along
the real analytic curve $\mathcal H^{+}_{\varepsilon}$ is bounded by  the cyclicity of the zero $h=0$ of the Abelian integral
$\tilde M_k$ .
Thus the imaginary part vanishes at most twice.
\end{itemize}
Summing up the above information we get at most five zeros.
 $\Box$

 \section{Proof of Theorem \ref{main}}
 \label{last}
 If the Bautin ideal $\mathbb B$ were principal, with generator  $\varepsilon = \varepsilon(\lambda)$, then we can write for the 
 displacement function 
$$
P_\lambda(h)-h =  \varepsilon c_1 (I_0(h) + O(\lambda)) +
\varepsilon c_2( I_2(h)+O(\lambda) +\varepsilon c_3    (I_4'(h) + O(\lambda) ) .
$$
which is the analogue of formula (\ref{returnmap}), and similar expressions hold true for the Dulac maps $d_\varepsilon^{\pm}$. Therefore we may repeat the arguments given in section \ref{sectionth4}, to produce exctly the same estimates for the zeros of the displacement map, as in the case of a one-parameter deformation. This would complete the proof of Theorem \ref{main}.

 Of course, he Bautin ideal is not principal, even if we localize it at the origin. Following \cite{gavr08}, we proceed to its principalization. Namely, consider the map 
\begin{equation}
\label{blowup}
 \R^4\to \R \mathbb P^2 : (\lambda_1,\lambda_2,\lambda_3,\lambda_4) \mapsto [\lambda_1:\lambda_2\lambda_3:\lambda_4]) 
\end{equation}
 which is well defined, except along the center variety $C= \{\lambda \in \C^4 : \lambda_1=\lambda_2\lambda_3=\lambda_4=0 \}$. By definition, the blow up $S$ of $\mathbb B$ is the Zarisky closure of the graph  of the map (\ref{blowup}). Clearly $S\subset \R^4\times \R \mathbb P^2$ is a singular algebraic surface of dimension four coming with natural projection (analytic map)
 $$
 \pi : S \to \R^4 .
 $$
The exceptional divisor of the blow up is the divisor $\pi^{-1}(C)$ which is a three-dimensional algebraic set having two irreducible components. Obviously $\pi^{-1}(0) = \mathbb P^2$ which is 
canonically identified to the the projectivized vector space of Melnikov functions $I_0, I_2, I_4'$, see \cite[section 3]{gafr}. The ideal $\mathbb B$ defines an ideal sheaf $\mathcal B$, and les $\pi^*\mathcal B$ be the 
 inverse image of $\mathcal B$ which is also an ideal sheaf this time on $S$. The main feature of the inverse ideal sheaf $\pi^*\mathcal B$ is that it is locally principal, see \cite[section 2]{gafr}. 
For instance, if $c=[c_1:c_2:c_3] \in \mathbb P^2$ with $c_1\neq 0$, then in a neighbourhood of $(0,[c_1:c_2:c_3] )$ on $S$ in which
$$
(\lambda, [\lambda_1:\lambda_2\lambda_3:\lambda_4]) \sim (0,[c_1:c_2:c_3]) 
$$
we may choose $\varepsilon = \lambda_1$ as a generator of the ideal and express $$\lambda_2\lambda_3 = \lambda_1 f_1(\lambda), \lambda_4 = \lambda_1 f_2(\lambda)$$
for suitable analytic $f_1,f_2$.

The above considerations show that for each point $(0,c)\in 0\times \pi^{-1}(0)\subset S $ we can find a neighbourhood in $S$ for which at most five limit cycles bifurcate from the 8-loop. 
This completes the proof of Theorem \ref{main}.

\end{document}